\long\def\symbolfootnote[#1]#2{\begingroup%
\def\thefootnote{\fnsymbol{footnote}}\footnote[#1]{#2}\endgroup}

\documentclass[12pt]{article}
\usepackage{amsfonts}

\usepackage{amsmath}

\newtheorem{theorem}{Theorem}

\begin{document}

\author{Ovidiu Munteanu}
\title{On a characterization of the complex hyperbolic space}
\date{February 1, 2008}
\maketitle

\begin{abstract}
Consider a compact K\"{a}hler manifold $M^m$ with Ricci curvature lower
bound $Ric_M\geq -2\left( m+1\right) .$ Assume that its universal cover $%
\widetilde{M}$ has maximal bottom of spectrum $\lambda _1\left( \widetilde{M}%
\right) =m^2.$ Then we prove that $\widetilde{M}$ is isometric to the
complex hyperbolic space $\Bbb{CH}^m.$
\end{abstract}

\section{Introduction}

\symbolfootnote[0]{Reaserch partially supported by NSF grant No. DMS-0503735}

Complete Riemannian manifolds with Ricci curvature lower bound have been the
object of study of many authors and there are very interesting results about
such manifolds. An important approach is to see how the spectrum of the
Laplacian interacts with the geometry of the manifold. A famous result that
we recall here is S.Y. Cheng's comparison theorem \cite{C}. If the Ricci
curvature of a complete noncompact Riemannian manifold $N^n$ of dimension $n$
is bounded from below by $Ric_N\geq -\left( n-1\right) ,$ then Cheng's
theorem asserts that the bottom of the spectrum of the Laplacian has an
upper bound $\lambda _1\left( N\right) \leq \frac{\left( n-1\right) ^2}%
4=\lambda _1\left( \Bbb{H}^n\right) .$ This result is sharp, but we should
point out that there are in fact many manifolds with maximal $\lambda _1,$
more examples can be found by considering hyperbolic manifolds $N=\Bbb{H}%
^n/\Gamma $ obtained by the quotient of $\Bbb{H}^n$ by a Kleinian group $%
\Gamma $ (\cite{S}). While in general we cannot determine the class of
manifolds with $\lambda _1$ achieving its maximal value, recently there has
been important progress in some directions.

P. Li and J. Wang have studied the structure at infinity of a complete
noncompact Riemannian manifold that has $Ric_N\geq -\left( n-1\right) $ and
maximal bottom of spectrum $\lambda _1\left( N\right) =\frac{\left(
n-1\right) ^2}4.$ They proved that either the manifold is connected at
infinity (i.e. it has one end) or it has two ends. In case it has two ends
then it must split as a warped product of a compact manifold with the real
line \cite{L-W2}. Their result has since been extended in many other
situations, e.g. K\"{a}hler manifolds, quaternionic K\"{a}hler manifolds or
locally symmetric spaces.

Recently X. Wang \cite{W} has obtained an interesting result in a different
setting. Suppose $N^n$ is a compact Riemannian manifold with $Ric_N\geq
-\left( n-1\right) .$ Consider $\pi :\widetilde{N}\rightarrow N$ its
universal cover and assume that $\lambda _1(\widetilde{N})=\frac{\left(
n-1\right) ^2}4.$ Then $\widetilde{N}$ is isometric to the hyperbolic space $%
\Bbb{H}^n.$

Wang proved this theorem using the notion of Kaimanovich entropy $\beta $, 
\[
\beta =-\lim_{t\rightarrow \infty }\frac 1t\int_{\widetilde{N}}p\left(
t,x,y\right) \log p\left( t,x,y\right) dy, 
\]
where $p$ denotes the heat kernel on $\widetilde{N}$, which also plays an
important role in our discussion.

It should be pointed out that in Wang's theorem if the manifold $N$ is assumed to have
negative curvature (and removing the lower bound on Ricci curvature assumption)
then stronger results are already known from the work of Ledrappier, Foulon,
Labourie, Besson, Courtois, Gallot \cite{L, F-L, B-C-G}. In this case it can
be proved that if $h$ denotes the volume entropy of $N$ defined by 
\[
h=\lim_{R\rightarrow \infty }\frac{\log Vol\left( B_p(R)\right) }R 
\]
and $\lambda _1\left( \widetilde{N}\right) =\frac 14h^2$ then $N$ is locally
symmetric.

However, Wang's theorem is quite powerful because it does not assume negative
curvature.

It is a natural question to investigate these issues on K\"{a}hler manifolds.

A first question that one should ask is if Cheng's estimate can be inproved
in this case. The model space that we work with is now the
complex hyperbolic space $\Bbb{CH}^m.$ Recently Li-Wang have proved \cite
{L-W1} that for a complete noncompact K\"{a}hler manifold $M^m$ of complex
dimension $m$ if the bisectional curvature is bounded from below by $%
BK_M\geq -1,$ then $\lambda _1\left( M\right) \leq m^2=\lambda _1\left( \Bbb{%
CH}^m\right) .$ They proved that in fact if the bottom of spectrum $\lambda
_1\left( M\right) $ achieves its maximal value, then the manifold is either
connected at infinity or it has two ends and in this latter case it is
diffeomorphic to the product of a compact manifold with the real line and
the K\"{a}hler metric on $M$ has a specialized form.

We recently improved (see \cite{M}) Li-Wang's results for complete K\"{a}hler
manifolds that have a Ricci curvature lower bound, $Ric_M\geq -2\left(
m+1\right) ,$ which is a weaker assumption than bisectional curvature lower
bound. To prove the estimate for $\lambda _1\left( M\right) $ and the
structure at infinity for manifolds with maximum $\lambda _1$ we used a new
argument, a sharp integral estimate for the gradient of a certain class of
harmonic functions. In this paper we will use our argument to estimate the
Kaimanovich entropy from above, which will imply the following result.

\begin{theorem}
Let $M^m$ be a compact K\"{a}hler manifold of complex dimension $m$ and with
Ricci curvature bounded from below by $Ric_M\geq -2\left( m+1\right) .$
Assume its universal cover $\pi :\widetilde{M}\rightarrow M$ has maximal
bottom of spectrum, $\lambda _1\left( \widetilde{M}\right) =m^2.$ Then $%
\widetilde{M}$ is isometric to the complex hyperbolic space $\Bbb{CH}^m.$
\end{theorem}

We want to comment now about the particular case when $M$ has negative
curvature.

For K\"{a}hler manifolds with bisectional curvature lower bound $BK_M\geq -1$
it follows from Li-Wang \cite{L-W1} that volume entropy verifies the sharp
estimate $h\leq 2m.$ So maximal bottom of spectrum in this case implies $%
\lambda _1=\frac 14h^2$.

However, for only Ricci curvature lower bound $Ric_M\geq -2\left( m+1\right) 
$ it is not known if $h\leq 2m,$ so it is not clear how to apply the
Besson-Courtois-Gallot theorem in the negative curvature case.

\bigskip \textbf{Aknowledgement}. The author would like to express his deep
gratitude to his advisor, Professor Peter Li, for continuous help,
encouragement and many valuable discussions.

\newpage

\section{Proof of the Theorem}

First, let us set the notation. We use the notations in \cite{L-W, M}. If 
$ds^2=h_{\alpha \bar{\beta}}dz^\alpha d\bar{z}^\beta $ is the K\"{a}hler
metric on $\widetilde{M},$ then $Re\left( ds^2\right) $ defines a Riemannian metric on $%
\widetilde{M}.$

Note that if $\left\{ e_1,e_{2},..,e_{2m}\right\} $ with $e_{2k}=Je_{2k-1}$
for $k\in \left\{ 1,2,..,m\right\} $ is an orthonormal frame with respect to
the Riemannian metric on $\widetilde{M}$ then $\left\{ v_1,..,v_m\right\} $ is a unitary
frame of $T_x^{1,0}\widetilde{M},$ where 
\[
v_k=\frac 12\left( e_{2k-1}-\sqrt{-1}e_{2k}\right) . 
\]
In this notation the following formulas hold 
\begin{eqnarray*}
\nabla f\cdot \nabla g &=&2\left( f_\alpha f_{\bar{\alpha}}+g_\alpha g_{\bar{%
\alpha}}\right) \\
\Delta f &=&4f_{\alpha \bar{\alpha}}.
\end{eqnarray*}
In the statement of the theorem, the Ricci curvature lower bound refers to
the Riemannian metric and it is equivalent to saying $Ric_{\alpha \bar{\beta}%
}\geq -\left( m+1\right) \delta _{\alpha \bar{\beta}}$ with respect to any
unitary frame.

To prove the theorem we follow the argument in \cite{W} and use the results in 
\cite{M}.

We first need to recall some facts about the Kaimanovich entropy.

There are a few equivalent formulations of this entropy. First, it can be
defined as a limit of the heat kernel: 
\[
\beta =\lim_{t\rightarrow \infty }\left( -\frac 1t\int_{\widetilde{M}%
}p\left( t,x,y\right) \log p\left( t,x,y\right) dy\right) , 
\]
where $p$ is the heat kernel on $\widetilde{M}.$ This definition is useful
because it can be showed that (a result of Ledrappier \cite{L}) 
\[
\beta \geq 4\lambda _1\left( \widetilde{M}\right) . 
\]
There is another very useful formula for $\beta ,$ using the minimal Martin
boundary of $\widetilde{M}.$ Let us quickly recall some known facts (see
e.g. \cite{A}).

Let $H\left( \widetilde{M}\right) $ denote the space of harmonic functions
on $\widetilde{M},$ with the topology of uniform convergence on compact
sets. Observe that $K_O=\left\{ u\in H\left( \widetilde{M}\right) :u\left(
O\right) =1,\;u>0\right\} $ is a compact and convex subset of $H\left( 
\widetilde{M}\right) $ so denote with $\partial ^{*}\widetilde{M}$ the set
of extremal points of $K_O$, i.e. points in $K_O$ that do not lie in any
open line segment in $K_O$. Note that a point of $K_O$ is extremal iff it is
a minimal harmonic function normalized at $O,$ therefore $\partial ^{*}%
\widetilde{M}$ is the minimal Martin boundary of $\widetilde{M}.$ Since $K_O$
is a metric space and it is compact and convex, by a theorem of Choquet it
results that for any positive harmonic function $h$ there is a unique Borel
measure $\mu ^h$ on the set of extremal points of $K_O$ so that 
\[
h\left( x\right) =\int_{\partial ^{*}\widetilde{M}}\xi \left( x\right) d\mu
^h\left( \xi \right) 
\]
In particular, for $h=1$ there exists a unique measure $\nu $ on $\partial
^{*}\widetilde{M}$ so that for any $x\in \widetilde{M},$ 
\[
\int_{\partial ^{*}\widetilde{M}}\xi \left( x\right) d\nu \left( \xi \right)
=1. 
\]
Let $\Gamma $ denote the group of deck transformations on $\widetilde{M},$
then there is a natural action of $\Gamma $ on $\partial ^{*}\widetilde{M}$,
defined by 
\[
\left( \gamma \xi \right) \left( x\right) =\frac{\xi \left( \gamma
^{-1}x\right) }{\xi \left( \gamma ^{-1}O\right) }, 
\]
for any $\xi \in \partial ^{*}\widetilde{M}$ and for any $\gamma \in \Gamma
. $

It is important to know how the measure $\nu $ is changed by the action of $%
\Gamma $ on $\partial ^{*}\widetilde{M}$, it can be easily seen that if $%
\eta =\gamma \xi ,$ then 
\[
\frac{d\nu \left( \eta \right) }{d\nu \left( \xi \right) }=\xi \left( \gamma
^{-1}O\right) . 
\]
For $x\in \widetilde{M}$ define 
\[
\omega \left( x\right) =\int_{\partial ^{*}\widetilde{M}}\xi ^{-1}\left(
x\right) \left| \nabla \xi \right| ^2\left( x\right) d\nu \left( \xi \right)
, 
\]
and notice that $\omega $ descends on $M.$ Indeed, for any $\gamma \in
\Gamma $ we have that 
\[
\left| \nabla \xi \right| ^2\left( \gamma x\right) =\left| \nabla \left(
\gamma ^{*}\xi \right) \right| ^2\left( x\right) , 
\]
where $\gamma ^{*}\xi $ is the pull back of $\xi ,$ i.e. $\gamma ^{*}\xi
=\xi \circ \gamma .$ Then it is easy to check using the Radon-Nikodym
derivative that for $\eta =\gamma ^{-1}\xi $ we have 
\[
\xi ^{-1}\left( \gamma x\right) \left| \nabla \xi \right| ^2\left( \gamma
x\right) d\nu \left( \xi \right) =\eta ^{-1}\left( x\right) \left| \nabla
\eta \right| ^2\left( x\right) d\nu \left( \eta \right) . 
\]
Then it clearly follows that 
\begin{eqnarray*}
\omega \left( \gamma x\right) &=&\int_{\partial ^{*}\widetilde{M}}\xi
^{-1}\left( \gamma x\right) \left| \nabla \xi \right| ^2\left( \gamma
x\right) d\nu \left( \xi \right) \\
&=&\int_{\partial ^{*}\widetilde{M}}\eta ^{-1}\left( x\right) \left| \nabla
\eta \right| ^2\left( x\right) d\nu \left( \eta \right) \\
&=&\omega \left( x\right) .
\end{eqnarray*}
We have showed that in fact $\omega $ is a well defined function on $M$.
This function can be used now to give another formula for the Kaimanovich
entropy. Everywhere in this paper we will denote by $dv$ the normalized
Riemannian volume form i.e. 
\[
dv=\frac 1{\int_M\sqrt{g}dx}\left( \sqrt{g}dx\right) . 
\]
By a formula of Kaimanovich (\cite{K}, also see \cite{L, W} ) the entropy can also be expressed as 
\begin{eqnarray*}
\beta &=&\int_M\omega dv \\
&=&\int_M\left( \int_{\partial ^{*}\widetilde{M}}\xi ^{-1}\left( x\right)
\left| \nabla \xi \right| ^2\left( x\right) d\nu \left( \xi \right) \right)
dv.
\end{eqnarray*}
So we have the following 
\begin{equation}
4\lambda _1\left( \widetilde{M}\right) \leq \int_M\left( \int_{\partial ^{*}%
\widetilde{M}}\xi ^{-1}\left( x\right) \left| \nabla \xi \right| ^2\left(
x\right) d\nu \left( \xi \right) \right) dv.  \label{1}
\end{equation}
For Riemannian manifolds X.\ Wang has used this inequality together with the
sharp Yau's gradient estimate (\cite{L-W2}) to prove his result in the Riemannian
setting.

For our problem, a sharp pointwise gradient estimate for K\"{a}hler
manifolds is not known to be true, but we know a way to obtain a sharp
integral estimate for the gradient of harmonic functions. So the goal is to
show that 
\[
\int_M\left( \int_{\partial ^{*}\widetilde{M}}\xi ^{-1}\left( x\right)
\left| \nabla \xi \right| ^2\left( x\right) d\nu \left( \xi \right) \right)
dv\leq 4m^2. 
\]
To show this, we use the argument in \cite{M}. The main technical point now
is to justify integration by parts (and in what sense) that was used in \cite
{M}.

Let $u=\log \xi ,$ then a simple computation shows that 
\[
u_{\alpha \overline{\beta }}=\xi ^{-1}\xi _{\alpha \overline{\beta }}-\xi
^{-2}\xi _\alpha \xi _{\overline{\beta }}. 
\]
For a fixed $x\in \widetilde{M}$ consider 
\[
\int_{\partial ^{*}\widetilde{M}}\xi \left( x\right) \left| u_{\alpha 
\overline{\beta }}\right| ^2\left( x\right) d\nu \left( \xi \right) . 
\]
We first claim that this integral is a finite number (depending on $x$).
Indeed, since $\partial ^{*}\widetilde{M}$ is compact and $d\nu $ is a
finite measure, it suffices to show the integrand is bounded. But this is
true because for fixed $x$ we can bound $\left| \xi _{\alpha \overline{\beta 
}}\right| \left( x\right) \leq C\left( x\right) \xi \left( O\right) =C\left(
x\right) .$ This can be seen as follows. Consider $B_O(R)$ a geodesic ball
of radius $R$ big enough so that $x\in B_O(R).$ Note that there exists a
constant $A>0$ so that $\Delta \left| \xi _{\alpha \overline{\beta }}\right|
\geq -A\left| \xi _{\alpha \overline{\beta }}\right| $ on $B_O(R).$ Such a
constant $A$ can be chosen to depend on the lower bound of the bisectional
curvature on $B_O(R),$ using the Bochner formula. Using now the mean value
inequality we get that there exists a constant $C_1$ depending on $R$ and $A$
so that 
\[
\left| \xi _{\alpha \overline{\beta }}\right| ^2\left( x\right) \leq
C_1\int_{B_O(R)}\left| \xi _{\alpha \overline{\beta }}\right| ^2. 
\]
It is known that by using integration by parts and suitable cut-off
functions that there exists a constant $C_2$ so that 
\[
\int_{B_O(R)}\left| \xi _{\alpha \overline{\beta }}\right| ^2\leq
C_2\int_{B_O(2R)}\xi ^2. 
\]
The right side of this inequality can now be bounded by $C_3\xi ^2\left(
O\right) ,$ using the Harnack inequality. Obviously, these constants will
depend on $R,$ nevertheless it follows that for $x$ fixed $\left| \xi
_{\alpha \overline{\beta }}\right| \left( x\right) $ will be bounded
uniformly for all $\xi ,$ which was our claim.

The second claim is that the function thus obtained actually descends on $M.$
This claim can be showed as above, now using the fact that since $M$ is
K\"{a}hler, the deck transformations are holomorphic, therefore for $\gamma
\in \Gamma $ and $\gamma ^{*}\xi $ the pull back of $\xi $ we have 
\[
\left| \left( \log \xi \right) _{\alpha \overline{\beta }}\right| ^2\left(
\gamma x\right) =\left| \left( \log (\gamma ^{*}\xi )\right) _{\alpha \bar{%
\beta}}\right| ^2\left( x\right) . 
\]
The rest of the proof follows the same line as for the gradient of $\xi $
(see above).

Therefore it makes sense to consider the following quantity: 
\begin{gather*}
\int_M\int_{\partial ^{*}\widetilde{M}}\xi \left( x\right) \left| u_{\alpha 
\overline{\beta }}\right| ^2\left( x\right) d\nu \left( \xi \right)
dv=\int_M\int_{\partial ^{*}\widetilde{M}}\xi ^{-1}\left( x\right) \left|
\xi _{\alpha \overline{\beta }}\right| ^2\left( x\right) d\nu \left( \xi
\right) dv \\
-2\int_M\int_{\partial ^{*}\widetilde{M}}\xi ^{-2}\left( x\right) (\xi
_{\alpha \overline{\beta }}\xi _{\overline{\alpha }}\xi _\beta )\left(
x\right) d\nu \left( \xi \right) dv \\
+\frac 1{16}\int_M\int_{\partial ^{*}\widetilde{M}}\xi ^{-3}\left( x\right)
\left| \nabla \xi \right| ^4\left( x\right) d\nu \left( \xi \right) dv,
\end{gather*}
where each of the integrals in the right side are also well defined by a
similar discussion.

We now want to justify integration by parts to show that 
\[
\int_M\int_{\partial ^{*}\widetilde{M}}\xi ^{-1}\left( x\right) \left| \xi
_{\alpha \overline{\beta }}\right| ^2\left( x\right) d\nu \left( \xi \right)
dv=\int_M\int_{\partial ^{*}\widetilde{M}}\xi ^{-2}\left( x\right) (\xi
_{\alpha \overline{\beta }}\xi _{\overline{\alpha }}\xi _\beta )\left(
x\right) d\nu \left( \xi \right) dv
\]
Consider $\left( U_i\right) $ a covering of $M$ with small open sets and let 
$\rho _i$ be a partition of unity subordinated to this covering. We can
choose $\left( U_i\right) $ so that each $U_i$ is diffeomorphic to an open
set $\widetilde{U}_i\subset \widetilde{M}$ via $\pi .$ We then have 
\begin{gather*}
\int_M\int_{\partial ^{*}\widetilde{M}}\xi ^{-1}\left( x\right) \left| \xi
_{\alpha \overline{\beta }}\right| ^2\left( x\right) d\nu \left( \xi \right)
dv \\
=\int_M\int_{\partial ^{*}\widetilde{M}}\xi ^{-1}\left( x\right) \xi
_{\alpha \overline{\beta }}\left( x\right) \left( \xi _{\bar{\alpha}}\left(
x\right) \sum_i\rho _i\left( \pi \left( x\right) \right) \right) _\beta d\nu
\left( \xi \right) dv \\
=\sum_i\int_M\int_{\partial ^{*}\widetilde{M}}\xi ^{-1}\left( x\right) \xi
_{\alpha \overline{\beta }}\left( x\right) \left( \xi _{\bar{\alpha}}\left(
x\right) \rho _i\left( \pi \left( x\right) \right) \right) _\beta d\nu
\left( \xi \right) dv \\
=\sum_i\int_{U_i}\int_{\partial ^{*}\widetilde{M}}\xi ^{-1}\left( x\right)
\xi _{\alpha \overline{\beta }}\left( x\right) \left( \xi _{\bar{\alpha}%
}\left( x\right) \rho _i\left( \pi \left( x\right) \right) \right) _\beta
d\nu \left( \xi \right) dv \\
=\sum_i\int_{\widetilde{U}_i}\int_{\partial ^{*}\widetilde{M}}\xi
^{-1}\left( x\right) \xi _{\alpha \overline{\beta }}\left( x\right) \left(
\xi _{\bar{\alpha}}\left( x\right) \rho _i\left( \pi \left( x\right) \right)
\right) _\beta d\nu \left( \xi \right) dv \\
=-\sum_i\int_{\widetilde{U}_i}\int_{\partial ^{*}\widetilde{M}}\left( \xi
^{-1}\left( x\right) \xi _{\alpha \overline{\beta }}\left( x\right) \right)
_\beta \left( \xi _{\bar{\alpha}}\left( x\right) \rho _i\left( \pi \left(
x\right) \right) \right) d\nu \left( \xi \right) dv \\
=\sum_i\int_{\widetilde{U}_i}\int_{\partial ^{*}\widetilde{M}}\xi
^{-2}\left( x\right) (\xi _{\alpha \overline{\beta }}\xi _{\overline{\alpha }%
}\xi _\beta )\left( x\right) \left( \xi _{\bar{\alpha}}\left( x\right) \rho
_i\left( \pi \left( x\right) \right) \right) d\nu \left( \xi \right) dv
\end{gather*}
\begin{eqnarray*}
&=&\sum_i\int_{U_i}\int_{\partial ^{*}\widetilde{M}}\xi ^{-2}\left( x\right)
(\xi _{\alpha \overline{\beta }}\xi _{\overline{\alpha }}\xi _\beta )\left(
x\right) \left( \xi _{\bar{\alpha}}\left( x\right) \rho _i\left( \pi \left(
x\right) \right) \right) d\nu \left( \xi \right) dv \\
&=&\sum_i\int_M\int_{\partial ^{*}\widetilde{M}}\xi ^{-2}\left( x\right)
(\xi _{\alpha \overline{\beta }}\xi _{\overline{\alpha }}\xi _\beta )\left(
x\right) \left( \xi _{\bar{\alpha}}\left( x\right) \rho _i\left( \pi \left(
x\right) \right) \right) d\nu \left( \xi \right) dv \\
&=&\int_M\int_{\partial ^{*}\widetilde{M}}\xi ^{-2}\left( x\right) (\xi
_{\alpha \overline{\beta }}\xi _{\overline{\alpha }}\xi _\beta )\left(
x\right) d\nu \left( \xi \right) dv.
\end{eqnarray*}
Let us mark out that everywhere in this formulas (and in the paper) a priori
the integrals on the minimal Martin boundary are taken for any (parameter) $%
x\in \widetilde{M}.$ Then, it can be justified that in fact these integrals
on $\partial ^{*}\widetilde{M}$ are invariant by the group of deck
transformations, so they are well defined functions on $M$. With this in
mind, in the third line above one should also justify that for each $i$ the
functions on $\widetilde{M}$ defined by $x\rightarrow \int_{\partial ^{*}%
\widetilde{M}}\xi ^{-1}\left( x\right) \xi _{\alpha \overline{\beta }}\left(
x\right) \left( \xi _{\bar{\alpha}}\left( x\right) \rho _i\left( \pi \left(
x\right) \right) \right) _\beta d\nu \left( \xi \right) $ descend on $M.$
This can be done by the same argument, and using that $\gamma ^{*}\left(
\rho _i\circ \pi \right) =$ $\rho _i\circ \pi ,$ for any $\gamma \in \Gamma .
$ It is also important that the function in $\xi $ which is integrated on
the minimal Martin boundary (for example $\xi \rightarrow \xi ^{-1}\xi
_{\alpha \bar{\beta}}\left( \xi _{\bar{\alpha}}\rho _i\circ \pi \right)
_\beta $) be homogeneous of degree 1 in $\xi .$ Thus we want to remark that
not quite any integration by parts is allowed by this procedure of lifting
the integrals on the universal covering.

This argument will be applied below every time we integrate by parts, it is
easy to check that the argument works in each case.

To simplify the writing, we will henceforth omit to write the argument $x$
and the measure $d\nu ,$ but we always assume the integrals on $\partial ^{*}%
\widetilde{M}$ are taken with respect to $d\nu $ and that all the functions
integrated on $\partial ^{*}\widetilde{M}$ depend on $x\in \widetilde{M}.$
For each of these integrals on the minimal Martin boundary it can be
justified that it is invariant by the group of deck transformations so it
legitimately defines a function on $M.$

We have thus proved that

\begin{eqnarray*}
\int_M\int_{\partial ^{*}\widetilde{M}}\xi \left| u_{\alpha \overline{\beta }%
}\right| ^2 &=&-\int_M\int_{\partial ^{*}\widetilde{M}}\xi ^{-2}(\xi
_{\alpha \overline{\beta }}\xi _{\overline{\alpha }}\xi _\beta ) \\
&&+\frac 1{16}\int_M\int_{\partial ^{*}\widetilde{M}}\xi ^{-3}\left| \nabla
\xi \right| ^4.
\end{eqnarray*}
Let us use again integration by parts to see that 
\begin{eqnarray}
-\int_M\int_{\partial ^{*}\widetilde{M}}\xi ^{-2}(\xi _{\alpha \overline{%
\beta }}\xi _{\overline{\alpha }}\xi _\beta ) &=&\int_M\int_{\partial ^{*}%
\widetilde{M}}\xi _\alpha \left( \xi ^{-2}\xi _{\overline{\alpha }}\xi
_\beta \right) _{\overline{\beta }}  \label{2} \\
&=&-\frac 18\int_M\int_{\partial ^{*}\widetilde{M}}\xi ^{-3}\left| \nabla
\xi \right| ^4+\int_M\int_{\partial ^{*}\widetilde{M}}\xi ^{-2}\xi _{%
\overline{\alpha }\overline{\beta }}\xi _\alpha \xi _\beta .  \nonumber
\end{eqnarray}
Note that the following inequality holds on $\widetilde{M}$: 
\[
\left| \xi _{\overline{\alpha }\overline{\beta }}\xi _\alpha \xi _\beta
\right| \leq \frac 14\left| \xi _{\alpha \beta }\right| \left| \nabla \xi
\right| ^2 
\]
so that we get 
\begin{gather}
2\int_M\int_{\partial ^{*}\widetilde{M}}\xi ^{-2}\xi _{\overline{\alpha }%
\overline{\beta }}\xi _\alpha \xi _\beta \leq \int_M\int_{\partial ^{*}%
\widetilde{M}}2\left( \xi ^{-1/2}\left| \xi _{\alpha \beta }\right| \right)
\left( \frac 14\xi ^{-3/2}\left| \nabla \xi \right| ^2\right)  \nonumber \\
\leq \frac m{m+1}\int_M\int_{\partial ^{*}\widetilde{M}}\xi ^{-1}\left| \xi
_{\alpha \beta }\right| ^2+\frac 1{16}\frac{m+1}m\int_M\int_{\partial ^{*}%
\widetilde{M}}\xi ^{-3}\left| \nabla \xi \right| ^4.  \label{3}
\end{gather}
Moreover, again integrating by parts we have 
\begin{eqnarray*}
\int_M\int_{\partial ^{*}\widetilde{M}}\xi ^{-1}\left| \xi _{\alpha \beta
}\right| ^2 &=&\int_M\int_{\partial ^{*}\widetilde{M}}\xi ^{-2}\xi _{%
\overline{\alpha }\overline{\beta }}\xi _\alpha \xi _\beta
-\int_M\int_{\partial ^{*}\widetilde{M}}\xi ^{-1}\xi _\alpha \xi _{\bar{%
\alpha}\bar{\beta}\beta } \\
&\leq &\int_M\int_{\partial ^{*}\widetilde{M}}\xi ^{-2}\xi _{\overline{%
\alpha }\overline{\beta }}\xi _\alpha \xi _\beta +\frac{m+1}%
4\int_M\int_{\partial ^{*}\widetilde{M}}\xi ^{-1}\left| \nabla \xi \right|
^2,
\end{eqnarray*}
using that $\xi $ is harmonic, the Ricci identities and the lower bound of
the Ricci curvature: 
\begin{eqnarray*}
-\xi _\alpha \xi _{\bar{\alpha}\bar{\beta}\beta } &=&-\xi _\alpha \xi _{\bar{%
\beta}\bar{\alpha}\beta } \\
&=&-\xi _\alpha \xi _{\bar{\beta}\beta \bar{\alpha}}-Ric_{\alpha \bar{\beta}%
}\xi _{\bar{\alpha}}\xi _\beta \\
&=&-Ric_{\alpha \bar{\beta}}\xi _{\bar{\alpha}}\xi _\beta \\
&\leq &\left( m+1\right) \xi _\alpha \xi _{\bar{\alpha}} \\
&=&\frac{m+1}4\left| \nabla \xi \right| ^2.
\end{eqnarray*}

Plug this inequality into (\ref{3}) and it follows 
\begin{gather*}
\frac{m+2}{m+1}\int_M\int_{\partial ^{*}\widetilde{M}}\xi ^{-2}\xi _{%
\overline{\alpha }\overline{\beta }}\xi _\alpha \xi _\beta \leq \frac
m4\int_M\int_{\partial ^{*}\widetilde{M}}\xi ^{-1}\left| \nabla \xi \right|
^2 \\
+\frac 1{16}\frac{m+1}m\int_M\int_{\partial ^{*}\widetilde{M}}\xi
^{-3}\left| \nabla \xi \right| ^4.
\end{gather*}
Getting back to (\ref{2}) we obtain 
\begin{gather*}
-\int_M\int_{\partial ^{*}\widetilde{M}}\xi ^{-2}\xi _{\alpha \overline{%
\beta }}\xi _{\overline{\alpha }}\xi _\beta \leq \left( -\frac 18+\frac 1{16}%
\frac{\left( m+1\right) ^2}{m\left( m+2\right) }\right) \int_M\int_{\partial
^{*}\widetilde{M}}\xi ^{-3}\left| \nabla \xi \right| ^4 \\
+\frac{m(m+1)}{4\left( m+2\right) }\int_M\int_{\partial ^{*}\widetilde{M}%
}\xi ^{-1}\left| \nabla \xi \right| ^2
\end{gather*}
We have thus proved that 
\begin{gather*}
\int_M\int_{\partial ^{*}\widetilde{M}}\xi \left| u_{\alpha \overline{\beta }%
}\right| ^2\leq \frac 1{16}\frac 1{m\left( m+2\right) }\int_M\int_{\partial
^{*}\widetilde{M}}\xi ^{-3}\left| \nabla \xi \right| ^4 \\
+\frac{m(m+1)}{4\left( m+2\right) }\int_M\int_{\partial ^{*}\widetilde{M}%
}\xi ^{-1}\left| \nabla \xi \right| ^2.
\end{gather*}
The estimate from below is straightforward: 
\[
\left| u_{\alpha \overline{\beta }}\right| ^2\geq \sum_\alpha \left|
u_{\alpha \bar{\alpha}}\right| ^2\geq \frac 1m\left| \sum_\alpha u_{\alpha 
\overline{\alpha }}\right| ^2=\frac 1{16m}\xi ^{-4}\left| \nabla \xi \right|
^4. 
\]
Hence, this shows that 
\begin{equation}
\int_M\int_{\partial ^{*}\widetilde{M}}\xi ^{-3}\left| \nabla \xi \right|
^4\leq 4m^2\int_M\int_{\partial ^{*}\widetilde{M}}\xi ^{-1}\left| \nabla \xi
\right| ^2.  \label{4}
\end{equation}
Finally, using the Schwarz inequality and the fact that $\int_{\partial ^{*}%
\widetilde{M}}\xi =1$ we get 
\begin{eqnarray*}
\int_M\int_{\partial ^{*}\widetilde{M}}\xi ^{-1}\left| \nabla \xi \right| ^2
&\leq &\left( \int_M\int_{\partial ^{*}\widetilde{M}}\xi ^{-3}\left| \nabla
\xi \right| ^4\right) ^{\frac 12}\left( \int_M\int_{\partial ^{*}\widetilde{M%
}}\xi \right) ^{\frac 12} \\
&=&\left( \int_M\int_{\partial ^{*}\widetilde{M}}\xi ^{-3}\left| \nabla \xi
\right| ^4\right) ^{\frac 12}.
\end{eqnarray*}
Combined with (\ref{4}) and (\ref{1}) this gives indeed that 
\[
4\lambda _1\left( \widetilde{M}\right) \leq \int_M\left( \int_{\partial ^{*}%
\widetilde{M}}\xi ^{-1}\left( x\right) \left| \nabla \xi \right| ^2\left(
x\right) d\nu \left( \xi \right) \right) dv\leq 4m^2, 
\]
as claimed.

Since we know $\lambda _1\left( \widetilde{M}\right) =m^2,$ it follows that
all inequalities used in this proof will be (pointwise) equalities on $%
\widetilde{M}$ for almost all $\xi \in \partial ^{*}\widetilde{M}.$ Indeed,
this is true because everywhere in our proof the inequalities were proved by
integrating on $\partial ^{*}\widetilde{M}$ some inequalities at $x\in 
\widetilde{M}$ that hold for each $\xi \in \partial ^{*}\widetilde{M}$.

Tracing back our argument, in \cite{M} we proved that for $B=\frac
1{2m}\log \xi $ we have 
\begin{eqnarray*}
\left| \nabla B\right| &=&1 \\
Hess_B\left( X,Y\right) &=&-g\left( X,Y\right) +g\left( \nabla B,X\right)
g(\nabla B,Y) \\
&&-g\left( J\nabla B,X\right) g\left( J\nabla B,Y\right)
\end{eqnarray*}
where $Hess_B$ denotes the real Hessian of $B$ .

From the work of Li-Wang \cite{L-W} we know that in this case, if the
manifold has bounded curvature then it is isometric to $\Bbb{CH}^m.$ This is
always the case for our setting, since $\widetilde{M}$ covers a compact
manifold, so its curvature is bounded. \textbf{Q.E.D.}

{\scriptsize DEPARTMENT OF MATHEMATICS, UNIVERSITY OF CALIFORNIA, IRVINE,
CA, 92697-3875}

{\small E-mail address: omuntean@math.uci.edu}

\end{document}